# SWASHES: A LIBRARY FOR BENCHMARKING IN HYDRAULICS

*SWASHES : Une bibliothèque de bancs d'essai en hydraulique*

## Olivier DELESTRE

Laboratoire J.A. Dieudonné & Polytech Nice – Sophia
Postal address[1]: UMR CNRS 7351 UNSA, 06108 Nice Cedex 02, France
e-mail : delestre@unice.fr

## Carine LUCAS

Laboratoire MAPMO – Fédération Denis Poisson
Postal address: UMR CNRS 7349, Université d'Orléans, B.P. 6759, 45067 Orléans Cedex 2, France
e-mail : carine.lucas@univ-orleans.fr

## Pierre-Antoine KSINANT

Institut National de la Recherche Agronomique (INRA)
Postal address: UR 0272 Science du sol, Centre de recherche d'Orléans, CS 40001, 45075 Orléans Cedex 2, France
e-mail : Pierre-Antoine.Ksinant@orleans.inra.fr

## Frédéric DARBOUX

Institut National de la Recherche Agronomique (INRA)
Postal address: UR 0272 Science du sol, Centre de recherche d'Orléans, CS 40001, 45075 Orléans Cedex 2, France
e-mail : Frederic.Darboux@orleans.inra.fr

## Christian LAGUERRE

Laboratoire MAPMO – Fédération Denis Poisson
Postal address: UMR CNRS 7349, Université d'Orléans, B.P. 6759, 45067 Orléans Cedex 2, France
e-mail : christian.laguerre@math.cnrs.fr

## François JAMES

Laboratoire MAPMO – Fédération Denis Poisson
Postal address: UMR CNRS 7349, Université d'Orléans, B.P. 6759, 45067 Orléans Cedex 2, France
e-mail : Francois.James@math.cnrs.fr

## Stéphane CORDIER

Laboratoire MAPMO – Fédération Denis Poisson
Postal address: UMR CNRS 7349, Université d'Orléans, B.P. 6759, 45067 Orléans Cedex 2, France
e-mail : stephane.cordier@math.cnrs.fr

De nombreux codes sont en cours d'élaboration pour résoudre les équations de Saint-Venant. Parce qu'elles sont utilisées dans les études hydrauliques et environnementales, leur capacité à simuler correctement les flux en eau est indispensable afin de préserver les infrastuctures et la sécurité humaine. Par conséquent, la validation de ces codes et des méthodes numériques associées est un problème essentiel. Des solutions analytiques de référence constitueraient une excellente réponse à ces questions. Toutefois, les solutions analytiques aux équations de Saint-Venant sont rares. Et surtout, elles ont été publiées sur une période de plus de cinquante ans, ce qui fait qu'elles sont dispersées à travers la littérature.
Dans cet article, un nombre important de solutions analytiques aux équations de Saint-Venant est décrit dans un formalisme unifié. Elles englobent une grande variété de conditions d'écoulement (supercritique, sous-critique, choc...), en 1 ou 2 dimensions d'espace, avec ou sans frottement, pluie et topographie, pour des écoulements transitoires ou à l'état stationnaire. Une caractéristique originale est que les codes source correspondants sont mis gratuitement à disposition de la communauté (http://www.univ-orleans.fr/mapmo/soft/SWASHES), afin que les utilisateurs de modèles en eaux peu profondes puissent facilement trouver un banc d'essai adaptable pour valider leurs méthodes numériques

---

[1] Corresponding author

*Numerous codes are being developed to solve Shallow Water equations. Because they are used in hydraulic and environmental studies, their capability to simulate properly flow dynamics is essential to guarantee infrastructure and human safety. Hence, validating these codes and the associated numerical methods is an important issue. Analytic solutions would be excellent benchmarks for these issues. However, analytic solutions to Shallow Water equations are rare. Moreover, they have been published on an individual basis over a period of more than five decades, making them scattered through the literature.*

*In this paper, a significant number of analytic solutions to the Shallow Water equations is described in a unified formalism. They encompass a wide variety of flow conditions (supercritical, subcritical, shock ...), in 1 or 2 space dimensions, with or without rain and soil friction, for transitory flow or steady state. An original feature is that the corresponding source codes are made freely available to the community (http://www.univ-orleans.fr/mapmo/soft/SWASHES), so that users of Shallow Water based models can easily find an adaptable benchmark library to validate their numerical methods.*

## Key words

## I   INTRODUCTION

Nowadays, Shallow-Water equations are widely used to model flows in various contexts, such as: overland flow [21,40], rivers [9,25], flooding [10,15], dam breaks [1,43], nearshore [6,33], tsunami [22,29,36]. This system of partial differential equations (PDE-s) proposed by Adhémar Barré de Saint-Venant in 1871 to model flows in a channel [4], consist in a system of conservation laws describing the evolution of the height and mean velocity of the fluid.

In real situations (realistic geometry, sharp spatial or temporal variations of the parameters), it is impossible to give an analytic *formulae* for the solutions of this system of PDE-s. Thus, there is a necessity to develop specific numerical schemes to compute approximate solutions for these PDE-s [7,30,42]. Implementation of such methods implies a subsequent step of code validation.

Validation of a model (that is the equations, the numerical methods and their implementation) is essential to know if it describes suitably the considered phenomena. At least three complementary kinds of numerical tests help us ensure that a numerical code is suitable for the considered system of equations. First, we can perform convergence or stability analysis (*e.g.* by refining the mesh). But this validates only the numerical method and its implementation. Second, approximate solutions can be compared with analytic solutions available for some simplified or specific cases. Finally, numerical results can be applied on experimental data, produced indoor or outdoor. This step should be done after the previous two; it is the most difficult one and must be validated by a specialist of the domain. In [18], we have focused on the second approach.

In numerical codes validation, analytic solutions seem to be underused. We think that there are two possible reasons. First, each analytic solution has a limited scope in terms of flow conditions. Second, as they are dispersed through the literature, they are difficult to find. However, an important number of published analytic solutions allow to embrace a wide range of flow conditions. Thus, overall, the existing analytic solutions have a large potential for numerical codes validation. In literature, we can find benchmarks for hydraulic river modelling software [14]. But they are too specific (weirs, pump, culverts... treatment validation) and there are too few tests available for numerical methods for Shallow-Water equations.

In [18], we have tried to overcome these problems both by gathering a significant set of analytic solutions and by providing the corresponding source codes. [18] describes the analytic solutions and gives some comments about their use and advantage. The source codes are freely available through the SWASHES (Shallow-Water Analytic Solutions for Hydraulic and Environmental Studies) library. With SWASHES software, we do not pretend to list all existing analytic solutions. Indeed, SWASHES is a framework to which users are invited to contribute by sending other analytic solutions together with the dedicated code.

The paper is organized as follows: in section II, we briefly present the notations we use and the main properties of Shallow-Water equations. In section III, we will focus on stationary solutions which are well-known by the hydraulics community but much less by mathematicians, *i.e.* "backwater curves". Lastly, in section IV, we will present SWASHES and the interest of solutions described in [18].

# II  EQUATIONS, NOTATIONS AND PROPERTIES

In the first section, we give the "complete" Shallow-Water system in two space dimensions, *i.e.* with topography, rain, infiltration, soil friction and viscous term. Then, we give this system in one space dimension and its main properties are recalled.

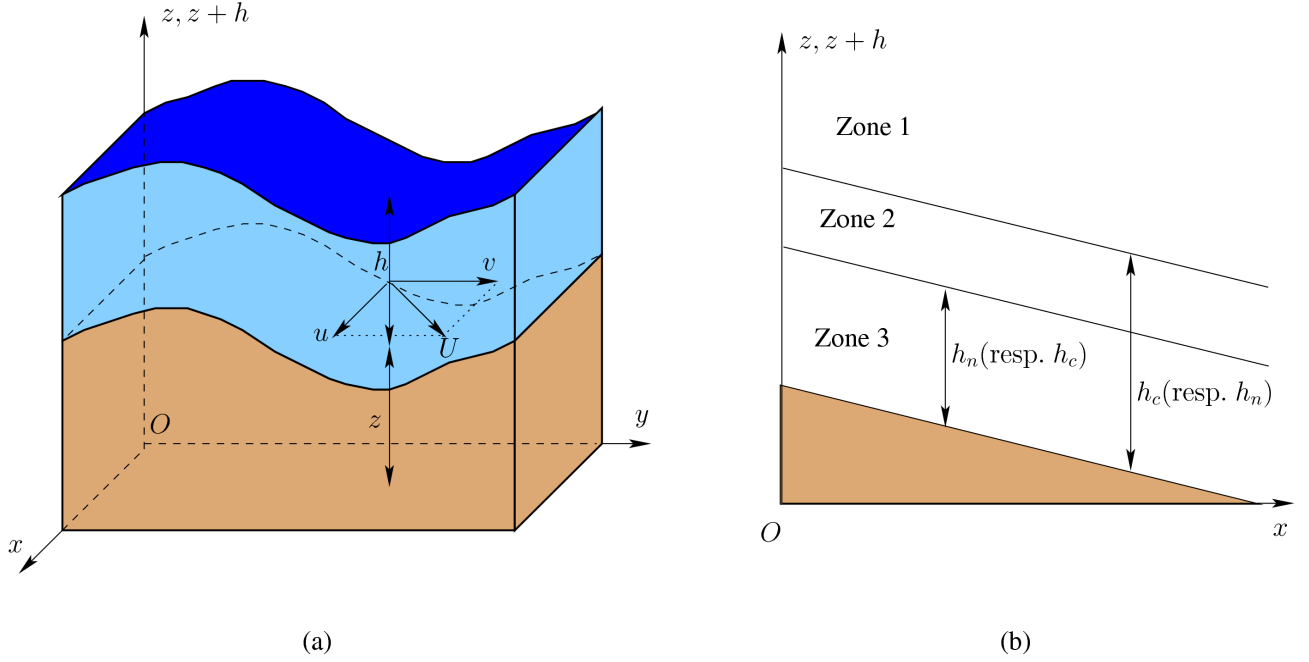

|                    (a)                    |                    (b)                    |

**Figure 1**: Illustration of variables of (a) Shallow Water equations and (b) zones for classification of free surface profiles.

## II.1   General settings

The Shallow-Water equations in two space dimensions take the following form:

$$\partial_t h + \partial_x (hu) + \partial_y (hv) = R - I$$
$$\partial_t (hu) + \partial_x \left( hu^2 + \frac{gh^2}{2} \right) + \partial_y (huv) = gh(S_{Ox} - S_{fx}) + \mu S_{dx} \quad , \qquad (1)$$
$$\partial_t (hv) + \partial_x (huv) + \partial_y \left( hu^2 + \frac{gh^2}{2} \right) = gh(S_{Oy} - S_{fy}) + \mu S_{dy}$$

where the unknowns of the model are the water height ($h(t,x,y)$ [L]) and $u(t,x,y)$, $v(t,x,y)$ the horizontal components of the vertically averaged velocity $[L/T]$ (**Figure 1 a**) and $g = 9.81 \, m/s^2$ is the gravity constant. The first equation is the mass conservation equation. The other two equations are the momentum conservation equations; they involve forces such as gravity and friction. We give now a short description of all the terms with their physical dimensions.

- $z$ is the topography $[L]$, since erosion is not considered here, it depends only on space, $z(x,y)$, and $S_{Ox}$ (resp. $S_{Oy}$) is the opposite of the slope in the $x$ (resp. $y$) direction, $S_{Ox} = -\partial_x z(x,y)$ (resp. $S_{Oy} = -\partial_y z(x,y)$);
- $R \geqslant 0$ is the rain intensity $[L/T]$, it is a given function $R(t,x,y) \geqslant 0$. In [17], it is considered uniform in space;
- $I$ is the infiltration rate $[L/T]$. It is given by another model (such as Horton, Philip, Green-Ampt, Richards …) and is not taken into account in the following;
- $S_f = (S_{fx}, S_{fy})$ is the friction force/law and may take various forms, depending on both soil and flow properties. In the *formulae* below, $U$ is the velocity vector $U = (u,v)$ with $|U| = \sqrt{u^2 + v^2}$ and $Q$ is the discharge $Q = (hu, hv)$. In hydrological models, two families of friction laws are encountered, based on empirical considerations. On one hand, we have the family of Manning-Strickler's friction laws

$$S_f = C_f \frac{U|U|}{h^{4/3}} = C_f \frac{Q|Q|}{h^{10/3}} \qquad (2)$$

$C_f = n^2$, where $n$ is the Manning's coefficient $[L^{-1/3}T]$.

On the other hand, the laws of Darcy-Weisbach's and Chézy's family write

$$S_f = C_f \frac{U|U|}{h} = C_f \frac{Q|Q|}{h^3} \quad . \qquad (3)$$

With $C_f = f/(8g)$, $f$ a dimensionless coefficient, (resp. $C_f = 1/C^2$, $C[L^{1/2}/T]$ ) we get the Darcy-Weisbach's (resp. Chézy's) friction law. Notice that the friction may depend on the space variable, especially for large areas. In the following this will not be considered.

- finally, $\mu S_d = (\mu S_{dx}, \mu S_{dy})$ is the viscous term with $\mu \geq 0$ the viscosity of the fluid $[L^2/T]$.

**II.2   Properties**

In this section, some properties of the Shallow-Water model are recalled. These properties are useful to the flow description. To simplify, we consider the one-dimensional model, but extensions to two dimensions are straightforward. The 2D Shallow-Water system (1) rewrites

$$\partial_t h + \partial_x(hu) = R - I$$
$$\partial_t(hu) + \partial_x\left(hu^2 + \frac{gh^2}{2}\right) = gh(S_{Ox} - S_f) + \mu \partial_x(h\partial_x u) \quad . \qquad (4)$$

The left-hand side of this system is the transport part of the model. It corresponds to the flow of an ideal fluid in a flat channel, without friction, rain or infiltration. It is in fact the model introduced by Saint-Venant in [4]. Several important properties of the flow are included in this model. The one-dimensional equations is rewritten using vectors form, in order to highlight these properties:

$$\partial_t W + \partial_x F(W) = 0, where\ W = \begin{pmatrix} h \\ hu \end{pmatrix}, F(W) = \begin{pmatrix} hu \\ hu^2 + \frac{gh^2}{2} \end{pmatrix} , \qquad (5)$$

with $F(W)$ the flux of the equation. With the following nonconservative form, where $A(W) = F'(W)$ is the jacobian matrix or the matrix of transport coefficients:

$$\partial_t W + A(W)\partial_x W = 0, with\ A(W) = F'(W) = \begin{pmatrix} u & 1 \\ u^2 + gh & 2u \end{pmatrix} , \qquad (6)$$

the transport is more clearly evidenced. More precisely, when, the matrix turns out to be diagonalizable, with eigenvalues

$$\lambda_1(W) = u - \sqrt{gh} < u + \sqrt{gh} = \lambda_2(W) \quad . \qquad (7)$$

In that case, the system is said to be strictly hyperbolic (see among others [23]). The eigenvalues (7) are the velocities of surface waves on the fluid. For dry zones (*i.e.* if $h = 0\,m$ ), the eigenvalues coincide. In that case, the system is no longer hyperbolic, and this induces difficulties at both theoretical and numerical levels, such as negative water depths... Designing numerical schemes that preserve positivity for $h$ is very important in this context.

With these *formulae,* we recover a useful classification of flows. Indeed if $|u| < \sqrt{gh}$, the characteristic velocities (of the fluid $u$ and of the waves $\sqrt{gh}$ ) have opposite signs, and informations propagate upward as well as downward the flow. The flow is said to be subcritical or fluvial. In the other case, when $|u| > \sqrt{gh}$, the flow is supercritical, or torrential, all the informations go downwards. A transcritical regime exists when some parts of a flow are subcritical, other supercritical.

Since we have two unknowns $h$ and $u$ (or equivalently $h$ and $q = hu$ ), a subcritical flow is therefore determined by one upstream and one downstream value, whereas a supercritical flow is completely determined by the two upstream values. Thus for numerical simulations, we have to impose one variable for subcritical inflow/outflow. For supercritical inflow, we impose both variables and free boundary conditions are considered (see for example [8,13,35]).

Two quantities allow us to determine the type of flow. The first one is a dimensionless parameter called the Froude number

$$Fr = \frac{|u|}{\sqrt{gh}} \quad . \qquad (8)$$

Its analogue in gas dynamics is the Mach number. If $Fr<1$ (resp. $Fr>1$), the flow is subcritical (resp. supercritical). The other essential quantity is the critical height $h_c$ which writes

$$h_c = \left(\frac{q}{\sqrt{g}}\right)^{2/3}, \qquad (9)$$

for a given discharge $q=hu$. It is a very readable criterion for criticality: the flow is subcritical (resp. supercritical) if $h>h_c$ (resp. $h<h_c$).

In presence of additional terms, we have to consider other properties, such as the occurrence of steady state (or equilibrium) solutions. In section III, we will focus on backwater curves which are specific steady state solutions.

## III BACKWATER CURVES

By considering system (4) at steady state ($\partial_t h = \partial_t u = \partial_t q$), without rain and diffusion ($R=0$ and $\mu=0$), it rewrites

$$q = q_0$$
$$\partial_x h = \frac{S_0 - S_f}{1 - Fr^2}. \qquad (10)$$

In hydraulics, equation (10) is used as a base for theoretical analysis of the water surface profiles obtained for different flow conditions in open channels [11,12,27]. It is called the gradually varied flow equation [11,12]. Water-surface profiles can be deduced theoretically and drawn qualitatively by studying the relative position of the profile ($h$) with respect to the critical-depth line ($h_c$) and the normal-depth line ($h_n$, the height solution of equation $S_0 - S_f = 0$). We notice that the normal height $h_n$ depends on the slope $S_0$ while the critical height $h_c$ does not depend on $S_0$. Bottom slopes are classified into five categories (designated by the first letter of the name): mild M if $h_n > h_c$, critical C if $h_n = h_c$, steep S if $h_n < h_c$, horizontal H if $S_0 = 0$ and adverse A if $S_0 < 0$. Now, we have to designate the relative position of the free-surface. In the cases of the mild and steep slopes, the space above the topography is divided into three regions by the normal height and the critical height (**Figure 1 b**). For the adverse, horizontal and critical slopes there are only two regions. Because the normal height does not exist for the two first slopes and is the same as the critical one for the critical slope. The region between the lower line and the topography is designated as "zone 3", the region between the upper and lower lines is designated as "zone 2" and that above both lines is designated as "zone 1". Thus, we have 13 different types of water height profiles. This technique allows to make qualitative observations about various types of free-surface profiles. These observations allow to draw the profile without any detailed calculations. For example, we know if the water height increases or decreases with distance, how the profiles end at downstream and upstream limits... It is possible to put end to end several backwater curves to obtain a complete profile. In some cases, it is possible to generate this method to section average model (*i.e.* to take into account the shape of the cross-section).

For engineering applications, it is necessary to compute the flow conditions. But the gradually flow equation (10) is nonlinear, and the dependence on $h$ is complicated, so analytical solution is not possible, so we have to use high order numerical methods [35]. We start the computations from downstream if the flow is subcritical and from upstream otherwise (see [11,12,27]). Some computer computer programs such as HEC-RAS are based on this method. In next section, we will describe SWASHES library.

## IV ADVANCED ANALYTIC SOLUTIONS

The Shallow-Water Analytic Solutions for Hydraulic and Environmental Software (SWASHES) is freely available to the community through the SWASHES repository hosted at http://www.univ-orleans.fr/mapmo/soft/SWASHES. It is distributed under CeCILL-V2 (GPL-compatible) free software licence. When running the software, user must specify the choice of the solution as well as the number of cells for the discretization of the chosen analytic solution. The solution is calculated and is redirected in a gnuplot-compatible ASCII file. SWASHES is written in object-oriented ISO C++ that allows to easily implement a new solution.

We claim that SWASHES can be a useful tool for developers of Shallow-Water codes to evaluate the performances and properties of their own code. Indeed, SWASHES has been created because we have been developing a free software for the resolution of Shallow-Water equations, namely FullSWOF [16,19], and we wanted to validate it against analytic solutions. With SWASHES, a wide range of flow conditions are available as developed in the following subsection: steady state solutions and transitory solutions.

## IV.1 Steady state solutions

In case of steady states ($\partial_t h = \partial_t u = \partial_t q = 0$), the one-dimensional Shallow-Water equations (4) reduce in the following system

$$q = Rx + q_0$$
$$\partial_x z = \frac{1}{gh}\left(\frac{q^2}{h^2} - gh\right)\partial_x h - S_f(h,q) + \frac{\mu}{gh}\partial_x\left(h\partial_x \frac{q}{h}\right) \quad . \quad (11)$$

System (11) enables to produce an infinity of analytic solutions. For these solutions, the strategy consists in choosing either a topography and getting the associated water height or a water height and deducing the associated topography. With the first approach, we can get for example the well-known solutions for flow over a bump [24,28], the backwater curves (section III)... With the second approach, we get the MacDonald's solutions (which are section averaged) [31,32] and all their one dimensional variants [16,17].

Since [5], it is well known that the topography source term treatment is a crucial point in preserving steady states. Thanks to the steady state solutions gathered in [18], one can check if the steady state at rest and dynamic steady state solutions are satisfied by the considered schemes. These solutions, integrated in SWASHES, cover a wide variety of flow conditions (fluvial, torrential, transcritical, with shock…). Moreover, different source terms (topography, friction, rain and diffusion) are taken into account which allow to validate each source terms treatment.

## IV.2 Transitory solutions

In previous subsection, we spoke about the steady-state solutions of SWASHES. These solutions can be used to check if the numerical methods are able to keep/catch steady-state flows. But even if the initial conditions differ from the expected steady-state, we do not have information about the transitory behaviour. Thus, transitory solutions are also included in SWASHES, such as the dam break solutions of increasing complexity [37,39,20], 1D and 2D Thacker's and variants solutions [41,34,38]. These solutions allow to test moving wet/dry transitions, moving shock, moving wet/dry transitions with friction...

## V CONCLUSIONS

We have developed SWASHES, a free tool for benchmarking in hydraulics. It is opened to user's contributions. We think that it might be useful for codes/numerical testing.

## VI ACKNOWLEGMENTS AND THANKS

This study is part of ANR METHODE project ANR-07-BLAN-0232 granted by the French National Agency for Research. The authors wish to thank V. Caleffi and A.-C. Boulanger for their collaboration.

## VII REFERENCES AND CITATIONS